\documentclass{article}
\usepackage[utf8]{inputenc}
\usepackage{amsmath,amssymb}
\usepackage{amsfonts}
\usepackage{algpseudocode}
\usepackage{algorithm}
\usepackage{enumitem}
\usepackage[margin=1in]{geometry}
\usepackage{natbib}
\usepackage{graphicx}
\usepackage{soul}
\usepackage{setspace}
\usepackage{listings}
\usepackage{color,soul}
\usepackage{bbm}
\usepackage{authblk}
\usepackage{subcaption}

\newcommand{\vect}[1]{\textbf{\textit{#1}}}

\newcommand{\E}{\mathrm{E}}
\newcommand{\Var}{\mathrm{Var}}
\newcommand{\Cov}{\mathrm{Cov}}

\DeclareMathOperator*{\argmin}{arg\,min}
\newtheorem{theorem}{Theorem}

\newtheorem{lemma}[theorem]{Lemma}

\newlist{condenum}{enumerate}{1} 
\setlist[condenum]{label=\bfseries Condition \arabic*., ref=\arabic*, wide}

\title{Exponential Consistency of M-estimators in Generalized Linear Mixed Models}

\author[1*]{Andrea Bratsberg}
\author[1]{Magne Thoresen}
\author[2]{Abhik Ghosh}
\affil[1]{Oslo Centre for Biostatistics and
Epidemiology, Department of Biostatistics,
University of Oslo}
\affil[2]{Indian Statistical Institute, Kolkata, India}
\affil[*]{\textit{email: a.m.bratsberg@medisin.uio.no}}

\begin{document}
\date{}

\maketitle

\begin{abstract}
Generalized linear mixed models are powerful tools for analyzing clustered data, 
where the unknown parameters are classically (and most commonly) estimated by the maximum likelihood and restricted maximum likelihood procedures.
However, since the likelihood based procedures are known to be highly sensitive to outliers, M-estimators have become popular as a means to obtain robust estimates under possible data contamination. 
In this paper, we prove that, for  sufficiently smooth general loss functions defining the M-estimators in generalized linear mixed models, 
the tail probability of the deviation between the estimated and the true regression coefficients  have an exponential bound. 
This implies an exponential rate of consistency of these M-estimators under appropriate assumptions, generalizing the existing exponential consistency results from univariate to multivariate responses. 
We have illustrated this theoretical result further for the special examples of the maximum likelihood estimator and the robust minimum density power divergence estimator, 
a popular example of model-based M-estimators, in the settings of linear and logistic mixed models, comparing it with the empirical rate of convergence through simulation studies.
\end{abstract}

\textbf{Keywords:} M-estimators, Generalized linear mixed models, Exponential consistency, Minimum density power divergence estimator.

\section{Introduction}

Linear and generalized linear mixed models (GLMMs) are powerful tools for analyzing a wide range of clustered data. Examples include repeated measurement data, longitudinal data, hierarchical and multilevel models, among others. The most common way of analyzing such data is by far the maximum likelihood and the restricted maximum likelihood estimation procedures. However, these are known to be highly sensitive to outliers, and a larger class of robust estimators, such as M-estimators which also contain the maximum likelihood estimators (MLE) as a special case, have become popular. All M-estimators 
are defined as the minimizer of some suitable loss function, e.g. the negative log-likelihood in the case of MLE. Existing robust M-estimators and their generalizations under the linear mixed models  include the robustified version of log-likelihood for
repeated measurement data by \cite{huggins1993robust}, the robust restricted maximum likelihood of \cite{richardson1995robust} and the SMDM-estimator of \cite{koller2013robust}. The popular minimum density
power divergence approach of \cite{basu1998robust} has also been extended to linear mixed models by \cite{saraceno2020robust} and is shown to either outperform or be competitive to other existing robust estimation procedures.
However, such robust M-estimators for GLMMs have received much less attention; some scattered works include the maximum quasi-likelihood and residual maximum quasi-likelihood estimator of \cite{yau2002robust} and the  computationally tractable robust maximum likelihood estimation procedure of \cite{sinha2004robust}.

In this paper, we prove that for a suitably defined class of M-estimators with sufficiently smooth loss functions, the tail probability of the deviation between the estimated and the true regression coefficients in GLMMs have an exponential bound, implying the exponential rate of consistency of these M-estimators under appropriate assumptions. This result is a generalization of a similar result in \cite{sis_in_generalized_linear_models}, with a more general scope as it includes multivariate responses from GLMMs, and paves the way for developments of such exponential consistency results for more general multivariate data.  We further illustrate our main exponential consistency results for popular M-estimators under the linear and the logistic mixed effects models. As we do not wish to propose any new (robust) estimators in this paper and only illustrate the applicability of our exponential consistency results, we verify the associated assumptions for two existing examples of M-estimators, namely the classical (non-robust) MLE and the most recently proposed robust minimum density power divergence estimator of \cite{saraceno2020robust}. Our theoretical result is further complemented by empirical simulations for these examples.

\section{M-estimators for generalized linear mixed models}
\label{sec:Mestforglmm}
We consider a generalized linear mixed regression model with multivariate response $\boldsymbol{y} = [y_1,\cdots,y_m]^T \in\mathcal{Y}\subseteq \mathbb{R}^m$
and a set of linear predictors $\boldsymbol{X}\boldsymbol{\beta}$, 
where $\boldsymbol{X}=[\boldsymbol{x}_1^T, \cdots, \boldsymbol{x}_m^T]^T$ for some explanatory variables $\boldsymbol{x}_j\in \mathcal{X}\subseteq\mathbb{R}^p$,
$j=1, \ldots, m$, and the regression coefficient $\boldsymbol{\beta}\in\mathbb{R}^p$. 
In most applications, $\boldsymbol{x}_j$s are i.i.d copies of some underlying $p$-dimensional explanatory variables; we assume that $\vect{X}$ is of full column rank. In addition, we consider a set of $q$ unobserved random factors $U_r$, $r = 1,...,q$, where, for simplicity and without loss of generality, we have implicitly set the number of factors for each random effect to 1. Let us further denote the corresponding design matrix of random effects $\vect{Z} \in \mathbb{R}^{m\times q}$, which is often a submatrix of the design matrix $\vect{X}$ (but may also contain some other covariates). 
The random effect $\vect{u} = [U_1,...,U_q] \in \mathbb{R}^{q}$ is assumed to follow a multivariate 
normal distribution with mean $\textbf{0}$ and a positive-definite covariance matrix $\vect{G}=\vect{G}(\boldsymbol{\eta})$, where $\boldsymbol{\eta} \in \mathbb{R}^{q(q+1)/2}$ is the (unknown) random-effect variance components. Let $\boldsymbol{\theta} = [\theta_1,\cdots,\theta_m] = \vect{X}\boldsymbol{\beta}+\vect{Z}\vect{u}$, so that $\theta_j = \vect{x}_j^T\boldsymbol{\beta}+\vect{z}_j^T\vect{u}$, where $\vect{z}_j^T$ are the rows of $\vect{Z}$. Given these fixed and random effects, we assume that the conditional distribution of the response $\vect{y}$ belongs to the exponential family of distributions, i.e.
\begin{align}
    f(y_j\lvert \vect{u},\boldsymbol{\beta},\phi) = \exp\big\{(y_j\theta_j-b(\theta_j))/a(\phi)+c(y_j,\phi)\big\}
    \label{EQ:exponentialfamily}
\end{align}
for known functions $a, b$ and $c$ and dispersion parameter $\phi$. We have omitted the explicit conditioning on $\vect{X}$ for notational simplicity and will assume that this is given throughout, unless otherwise stated. For simplicity, we will assume $\phi = 1$. Then, the mean of $\vect{y}$ is given by
\begin{align}
\E(y_j\lvert\vect{u}) =h^{-1}(\theta_j) = b'(\theta_j) 
\label{EQ:conditionalmeanofy}
\end{align}
where $h(\cdot)$ is some known link function. 
The variance of $y_j$ given $\vect{u}$ is 
\begin{align*}
    \text{Var}(y_j\lvert\vect{u})= \nu(\E(y_j\lvert\vect{u})) = b''(\theta_j),
\end{align*}
where $\nu(\cdot)$ is a known function. Since the the elements of $\vect{y}$ are independent if $\vect{u}$ is given, the marginal distribution of the response vector $\vect{y}$ is given by
\begin{align}
f(\vect{y};\boldsymbol{\beta},\boldsymbol{\eta}) = \int f(\vect{y}\lvert \vect{u},\boldsymbol{\beta})f_u(\vect{u}\lvert\boldsymbol{\eta})d\vect{u}= \int \prod_{j=1}^mf(y_j\lvert \vect{u},\boldsymbol{\beta})f_u(\vect{u}\lvert\boldsymbol{\eta})d\vect{u},
\label{EQ:marginalyi}
\end{align}
where $f_u(\vect{u}\lvert\boldsymbol{\eta})$ is the Gaussian probability density function of $\vect{u}$ given by
\begin{align}
    f_u(\vect{u}\lvert \boldsymbol{\eta}) = \frac{1}{(2\pi)^{q/2}\lvert \vect{G}(\boldsymbol{\eta})\rvert^{1/2}}\exp\big\{-\frac{1}{2}\vect{u}^T\vect{G}(\boldsymbol{\eta})^{-1}\vect{u}\big\}.
\label{EQ:mixedeffectsdistribution}
\end{align}
The expression in \eqref{EQ:marginalyi} can often not be written in a closed form and needs to be computed numerically. Now, suppose that we have observations $(\boldsymbol{y}_i, \boldsymbol{X}_i)$, $i=1, \ldots, n$, as independent realizations of $(\boldsymbol{y}, \boldsymbol{X})$, so that $\vect{y}_i = [y_{i1},\cdots,y_{im}]^T$.
We want to estimate the unknown parameters $ (\boldsymbol{\beta}, \boldsymbol{\eta})$ based on these $n$ multivariate samples, yielding a total of $N=nm$ observed data points.
We assume that the total effective sample size $N$ is larger than $p+q(q+1)/2$, the total number of unknown parameters.  
Then, a general class of M-estimators may be defined as the minimizer of the empirical average of a sufficiently smooth loss function
$\rho(\cdot)$. More explicitly,
for a given $\rho(\cdot)$ function, the corresponding M-estimator of $(\boldsymbol{\beta}, \boldsymbol{\eta})$ is given by 
\begin{eqnarray}
(\widehat{\boldsymbol{\beta}}, \widehat{\boldsymbol{\eta}})
= \argmin_{\boldsymbol{\beta}, \boldsymbol{\eta}} \mathbb{P}_n\rho(\vect{y},\vect{X}\boldsymbol{\beta}, \boldsymbol{\eta}),
\label{EQ:M-est}
\end{eqnarray}
where $\mathbb{P}_n$ denotes the empirical average based on $n$ observations, i.e. 
$\mathbb{P}_n\left[h(\boldsymbol{y}, \boldsymbol{X})\right] = n^{-1}\sum_{i=1}^n h(\boldsymbol{y}_i, \boldsymbol{X}_i)$ for any function $h$. 
Note that the above definition of the M-estimator is slightly more stringent than the estimating equations based definition of M-estimators,
but cover many important estimators. Let us define the corresponding best fitting population parameter as
\begin{eqnarray}(\boldsymbol{\beta}_0, \boldsymbol{\eta}_0) 
	= \argmin_{\boldsymbol{\beta}, \boldsymbol{\eta}} \E\rho(\vect{y},\vect{X}\boldsymbol{\beta}, \boldsymbol{\eta}),
\label{EQ:M-func}
\end{eqnarray}
where the expectation is taken with respect to both $\vect{y}$ and $\vect{X}$.
The $\rho$ function must be suitably chosen so that this best fitting parameter should coincide with the true parameter value when the true data generating distribution belongs to the model family. 
It is well known that the M-estimator $(\widehat{\boldsymbol{\beta}}, \widehat{\boldsymbol{\eta}})
$ is consistent for $(\boldsymbol{\beta}_0, \boldsymbol{\eta}_0)$ as $n\rightarrow\infty$ and $m$ remains fixed \citep[Chapter~6.2]{huber2011robust}.
In this paper, we will show exponential consistency for $\widehat{\boldsymbol{\beta}}$ under suitable conditions, assuming $m$ can also vary. The dimension of the parameter vector is assumed to be fixed, and the variance parameter $\boldsymbol{\eta}$ is considered to be a nuisance parameter, estimated consistently by the same M-estimation approach. 

\subsubsection*{Example 1: Maximum likelihood estimator}
The most common example of a (non-robust) M-estimator is the classical maximum likelihood estimator. When the random effects $\vect{u}$ are given, the observations from the same group are independent and the joint likelihood for $\vect{y}$ is given by 
\begin{align*}
L(\boldsymbol{\beta},\boldsymbol{\eta}\lvert \vect{y},\vect{X}) = \int \prod_{i=1}^mf(y_i\lvert\vect{u},\boldsymbol{\beta})f_u(\vect{u}\lvert \boldsymbol{\eta})d\vect{u},
\end{align*}
where  $f(y_i\lvert \vect{u},\boldsymbol{\beta})$ is given by \eqref{EQ:exponentialfamily} and $f_u(\vect{u}\lvert \boldsymbol{\eta})$ is given in \eqref{EQ:mixedeffectsdistribution}.
The MLE estimates $(\widehat{\boldsymbol{\beta}}_{\text{ML}}, \widehat{\boldsymbol{\eta}}_{\text{ML}})$ of $(\boldsymbol{\beta},\boldsymbol{\eta})$ are defined as the maximizer of the above likelihood, which is equivalent to minimizing the negative log-likelihood; formally we may thus define
\begin{align*}
(\widehat{\boldsymbol{\beta}}_{\text{ML}}, \widehat{\boldsymbol{\eta}}_{\text{ML}})
 =-\argmin_{\boldsymbol{\beta}, \boldsymbol{\eta}} \mathbb{P}_n \ln L(\boldsymbol{\beta},\boldsymbol{\eta}\lvert \vect{y}, \vect{X}),
\end{align*}
which is clearly an M-estimator corresponding to the loss function $\rho(\vect{y}, \vect{X}\boldsymbol{\beta},\boldsymbol{\eta}) = - \ln L(\boldsymbol{\beta},\boldsymbol{\eta}\lvert \vect{y}, \vect{X})$.
The estimating equations for the MLE then becomes 
\begin{align}
    \mathbb{P}_n \Big[\frac{\partial  \rho(\vect{y},\vect{X}\boldsymbol{\beta}, \boldsymbol{\eta})}{\partial \boldsymbol{\beta}}\Big] = 0 \quad \text{ and } \quad
    \mathbb{P}_n \Big[\frac{\partial  \rho(\vect{y},\vect{X}\boldsymbol{\beta}, \boldsymbol{\eta})}{\partial \boldsymbol{\eta}}\Big] = 0.
    \label{EQ:estimatingequationsMLE}
\end{align}
No closed-form solution exists for these estimating equations in general, so they need to be solved numerically by, e.g., Monte Carlo Newton-Raphson methods.

\subsubsection*{Example 2: Minimum density power divergence estimator}

The minimum density power divergence estimator (MDPDE) of \cite{basu1998robust} has been defined for linear mixed models recently by \cite{saraceno2020robust}, but has not yet been developed for GLMMs. \cite{basu1998robust} defined the density power divergence measure (DPD) $d_\alpha(g,f)$ between two probability densities $g$ and $f$, as 
\begin{align*}
    d_\alpha(g,f) &= \int \Big\{f^{1+\alpha} -\Big(1+\frac{1}{\alpha}\Big)f^\alpha g+\frac{1}{\alpha}g^{1+\alpha}\Big\}, \quad \text{ if } \alpha > 0,\\
    d_0(g,f) &= \int g\log\Big(\frac{g}{f}\Big) \quad \text{ if } \alpha = 0.
\end{align*}
In the setting of estimation, we let $g$ denote the density function of the true data generating distribution $G$, and we model $g$ by the parametric family of densities $\mathcal{F}_{\boldsymbol{\vartheta}} = \{f_{\boldsymbol{\vartheta}}: \boldsymbol{\vartheta} \in \Theta \}$. \cite{ghosh2013robust} extended the MDPDE to the case of independent non-homogeneous observations, where $\vect{y}_1,...,\vect{y}_n$ are independent, but each $\vect{y}_i \sim g_i$ with $g_1,...,g_n$ being potentially different densities with respect to some common dominating measure. This is the case of our GLMMs as well, where we model each individual $g_i$ by the model density of the form $f(\vect{y}; \boldsymbol{\beta},\boldsymbol{\eta})$ given in \eqref{EQ:marginalyi}. Then, following \cite{ghosh2013robust}, we can define the MDPDE of the parameters $\boldsymbol{\vartheta} = (\boldsymbol{\beta},\boldsymbol{\eta})$ under the GLMMs as the minimizer of the objective function
\begin{align}
  \frac{1}{n}\sum_{i=1}^n\Big[\int f_i(\vect{y}_i;\boldsymbol{\beta},\boldsymbol{\eta})^{1+\alpha}dy-\Big(1+\frac{1}{\alpha}\Big)f_i(\vect{y}_i;\boldsymbol{\beta},\boldsymbol{\eta})^\alpha\Big].
\label{EQ:MDPDEobjfunc}
\end{align}
This is again an M-estimator with model dependent loss function given by 
\begin{align}
\rho(\vect{y},\vect{X}\boldsymbol{\beta},\boldsymbol{\eta})= \int f(\vect{y};\boldsymbol{\beta},\boldsymbol{\eta})^{1+\alpha}dy-\Big(1+\frac{1}{\alpha}\Big)f(\vect{y};\boldsymbol{\beta},\boldsymbol{\eta})^\alpha.
\label{EQ:mdpdelossfunctionGeneral}
\end{align}
As with the MLE, the MDPDE estimating equations will be of the form \eqref{EQ:estimatingequationsMLE},  but with the above-mentioned $\rho$-function, which need to be solved numerically. \cite{saraceno2020robust} have studied many important properties of this MDPDE, including its consistency, asymptotic normality and robustness, under the linear mixed-effect models, and empirically illustrated its superiority over the other existing robust estimators. Such nice properties of the MDPDEs are expected to hold under any GLMMs as well, but their formal investigation would be important future works (which is out of the scope and aim of the present manuscript). In this paper, we will only examine the rate of consistency of the MDPDEs under the linear and logistic mixed models in order to support our theoretical results on exponential consistency of the MDPDE as an example of the class of M-estimators.

\section{Main Results: Exponential Consistency of the Regression Coefficients}
\label{regularityConditions}

In order to prove our main results on the exponential consistency of the M-estimator of the regression coefficients $\boldsymbol{\beta}$, we will start by specifying some notations and required regularity conditions under the set-up of GLMMs. In the following, we denote the gradient with respect to $(\boldsymbol{\beta},\boldsymbol{\eta})$ by $\nabla$
and the $\ell_2$ and $\ell_\infty$ norms by $||\cdot||$ and $||\cdot||_\infty$, respectively. 
Also, let $\boldsymbol{1}_m$ denote the $m$-vector with all entries equal to one. 
For any set $S$ of points $(\boldsymbol{y}, \boldsymbol{X})$, we denote the indicator function of this set by 
$\mathbbm{1}_{S}(\vect{y},\vect{X})$ which takes the value one if $(\boldsymbol{y}, \boldsymbol{X})\in S$ and zero otherwise. We now consider the following set of assumptions.

\begin{itemize}
	\item[\textit{A.1}]
Suppose that, for a suitable constant $B>0$, $\boldsymbol{\beta}$ is an interior point of a compact and convex set $ \boldsymbol{\mathcal{B}} = \{\boldsymbol{\beta}: \lvert\lvert \boldsymbol{\beta} - \boldsymbol{\beta}_0 \rvert\rvert< B\}.$

\item[\textit{A.2}] There exists a positive constant $C$, independent of $m$, such that 
$\E ||\boldsymbol{1}_m^T\boldsymbol{X}||^2 \leq m C$.

\item[\textit{A.3}] The loss function $\rho(\boldsymbol{y}, \boldsymbol{X}\boldsymbol{\beta},\boldsymbol{\eta})$ is convex in $\boldsymbol{\beta}$ for all  $\boldsymbol\beta \in \boldsymbol{\mathcal{B}}$. In addition, the generalized information matrix
    \begin{align*}
    \boldsymbol{I}(\boldsymbol{\beta}, \boldsymbol{\eta}) 
    = \E\left[\left(\nabla\rho(\boldsymbol{y}, \boldsymbol{X}\boldsymbol{\beta},\boldsymbol{\eta})\right)
    \left(\nabla \rho(\boldsymbol{y}, \boldsymbol{X}\boldsymbol{\beta},\boldsymbol{\eta})\right)^T\right]  
    \end{align*} 
    is finite and positive definite at  $(\boldsymbol{\beta}_0,\boldsymbol{\eta}_0)$.  Moreover, 
    $\|\boldsymbol{I}(\boldsymbol{\beta}, \boldsymbol{\eta}) \|$ is bounded from above.
\item[\textit{A.4}] There exists a positive constant $V$ such that, for all $\boldsymbol\beta \in \boldsymbol{\mathcal{B}} $, we have 
\begin{align*}
	\E\left[\rho(\boldsymbol{y}, \boldsymbol{X}\boldsymbol{\beta},\boldsymbol{\eta_0})-\rho(\boldsymbol{y}, \boldsymbol{X}\boldsymbol{\beta}_0,\boldsymbol{\eta_0})\right] \geq 
	V\|\boldsymbol{\beta} - \boldsymbol{\beta}_0\|^2.
\end{align*}

\item[\textit{A.5}] The function $\rho(\boldsymbol{y}, \boldsymbol{X}\boldsymbol{\beta},\boldsymbol{\eta})$ satisfies the Lipschitz condition in $\vect{X}\boldsymbol{\beta}$ over $\boldsymbol{\beta} \in \boldsymbol{\mathcal{B}}$ at $\boldsymbol{\eta}=\boldsymbol{\eta}_0$, i.e. for a constant $k_n>0$ and for any $\boldsymbol{\beta}, \boldsymbol{\beta}'\in \boldsymbol{\mathcal{B}}$, we have
\begin{align*}
	\lvert \rho(\boldsymbol{y}, \boldsymbol{X}\boldsymbol{\beta},\boldsymbol{\eta}_0)-\rho(\boldsymbol{y}, \boldsymbol{X}\boldsymbol{\beta}',\boldsymbol{\eta}_0)\rvert \mathbbm{1}_{\Lambda_{n}}(\vect{y},\vect{X})
	\leq k_n \lvert \boldsymbol{1}_m^T\boldsymbol{X}\boldsymbol{\beta} - \boldsymbol{1}_m^T\boldsymbol{X}\boldsymbol{\beta}' 
	\rvert\mathbbm{1}_{\Lambda_{n}}(\vect{y},\vect{X}),
\end{align*}
where $\Lambda_{n} = \left\{(\boldsymbol{y}, \boldsymbol{X}) : \|\boldsymbol{y}\|_\infty \leq K_n^*, \|\boldsymbol{X}\|_\infty \leq K_n  \right\}$
for sufficiently large constants $K_n, K_n^* >0$, and $K_n \to \infty$ as $n \to \infty$.

\item[\textit{A.6}] There exists a sufficiently large constant $D > 0$ such that, with $b = Dk_n(mC/n)^{1/2}/V$, we have
\begin{align*}
\sup_{\substack{\boldsymbol{\beta}\in\boldsymbol{\mathcal{B}},\\
			\|\boldsymbol{\beta} - \boldsymbol{\beta}_0\| \leq b}}
	\lvert \E[\rho(\boldsymbol{y}, \boldsymbol{X}\boldsymbol{\beta},\boldsymbol{\eta}_0)-\rho(\boldsymbol{y}, \boldsymbol{X}\boldsymbol{\beta}_0,\boldsymbol{\eta}_0)]
	(1-\mathbbm{1}_{\Lambda_{n}}(\vect{y},\vect{X})) \rvert \leq o(mC/n),
\end{align*}
where the constants  $C$, $V$, $k_n$ and $\Lambda_n$ are as defined previously in Assumptions A.2--A.5. 
\end{itemize}

Note that A.1 is a technical assumption and Assumption A.2 specifies that the covariates remains bounded, in expectation, at a certain rate; these are clearly quite common assumptions in literature although the rate of boundedness in A.2 is somewhat specific to our case.  Furthermore, note that, $\boldsymbol{I}(\boldsymbol{\beta}, \boldsymbol{\eta})$ in A.3 is the usual Fisher information matrix when the $\rho$-function is the negative log-likelihood in case of the MLE. So, Assumption A.3 is just the generalization of the standard assumption related to the information matrix for the general loss function $\rho(\cdot)$ of the M-estimators. Next, the convexity of $\rho(\cdot)$ in $\boldsymbol{\beta}$ ensures the existence of a best-fitting parameter value $\boldsymbol{\beta}_0$, which our M-estimator of $\boldsymbol{\beta}$ would converge to asymptotically with probability tending to one at an exponential rate. Assumption A.4 requires the loss function to be strongly convex around the best-fitting parameter $(\boldsymbol{\beta}_0,\boldsymbol{\eta}_0)$, and by doing so ensuring that there exists a unique local minimum in the neighbourhood $\boldsymbol{\mathcal{B}}$. It is thus overlapping with A.3 but not the same. Assumption A.5 and A.6 are appropriate extensions of Condition B in \cite{sis_in_generalized_linear_models}.

\begin{theorem}[Exponential tail bound]\label{THM:1}
Suppose that Assumptions A.1--A.6 hold under the set-up of generalized linear mixed models and a given $\rho$-function defining the M-estimator of the parameters as in \eqref{EQ:M-est}. Then, for any $t >0$ and large enough $n$, we have
\begin{align*}
    \text{Pr}\left(\sqrt{\frac{n}{m}}\left\|\widehat{\boldsymbol{\beta}}-\boldsymbol{\beta}_0\right\|\geq \frac{16k_nC^{1/2}}{V}(1+t)\right) 
    \leq \exp\left(-\frac{2Ct^2}{K_n^2}\right)+ n\text{Pr}\left(\Lambda_n^c\right)+\text{Pr}(\Omega_n^c),
\end{align*}
where $\Lambda_n$ is defined in Condition A.5 and $\Omega_n = \{(\vect{y},\vect{X}):  \mathbbm{P}_n\rho(\boldsymbol{y}, \boldsymbol{X}\widehat{\boldsymbol{\beta}},\boldsymbol{\eta}_0) \leq  \mathbbm{P}_n\rho(\boldsymbol{y}, \boldsymbol{X}{\boldsymbol{\beta}_0},\boldsymbol{\eta}_0) \}.$
\end{theorem}
\textit{Proof:}

We refer to Lemmas 2--4 in \citet{sis_in_generalized_linear_models} for the proof.
Let us define
\begin{align*}
    \mathbb{G}_1(B) = \sup_{\boldsymbol{\beta}\in \mathcal{B}(B)}\lvert (\mathbb{P}_{n}-E)\{\rho(\boldsymbol{y}, \boldsymbol{X}\boldsymbol{\beta},{\boldsymbol{\eta}}_0)-\rho(\boldsymbol{y}, \boldsymbol{X}\boldsymbol{\beta}_0,{\boldsymbol{\eta}}_0))\}\mathbbm{1}_{\Lambda_{n}}(\vect{y},\vect{X})\rvert,
\end{align*}

\noindent\textit{Step 1:}\\
By an application of Lemma 2 of Fan and Song (2010) (Symmetrization theorem), we get
\begin{align*}
    \E(\mathbb{G}_1(B)) \leq 2\E\Big[\sup_{\boldsymbol{\beta}\in \mathcal{B}(B)}\left\lvert \mathbb{P}_{n}\varepsilon\{\rho(\boldsymbol{y}, \boldsymbol{X}\boldsymbol{\beta},{\boldsymbol{\eta}}_0)-\rho(\boldsymbol{y}, \boldsymbol{X}\boldsymbol{\beta}_0,{\boldsymbol{\eta}}_0)\}\mathbbm{1}_{\Lambda_{n}}(\vect{y},\vect{X})\right\rvert\Big],
\end{align*}
where $\varepsilon =(\varepsilon_1,...,\varepsilon_{n})$ is a Rademacher sequence. 
By Lemma 3 of \citet{sis_in_generalized_linear_models} and Condition A.5, we have that $\E(\mathbb{G}_1(B))$ is further bounded above by
\begin{align*}
    4k_n\E\Big\{\sup_{\boldsymbol{\beta}\in \mathcal{B}(B)}\left\lvert 
    \mathbb{P}_{n}\varepsilon\boldsymbol{1}_m^T\vect{X}(\boldsymbol{\beta}-\boldsymbol{\beta}_0)\mathbbm{1}_{\Lambda_{n}}(\vect{y},\vect{X})
    \right\rvert \Big\}
\end{align*}
Next, by the Cauchy-Schwartz inequality, we get that
\begin{eqnarray}
    \E\Big\{\sup_{\boldsymbol{\beta}\in \mathcal{B}(B)}\left\lvert 
    \mathbb{P}_{n}\varepsilon\boldsymbol{1}_m^T\vect{X}(\boldsymbol{\beta}-\boldsymbol{\beta}_0)\mathbbm{1}_{\Lambda_{n}}(\vect{y},\vect{X})
    \right\rvert \Big\} 
    &\leq& \E\|\mathbb{P}_{n}\varepsilon\vect{X}\boldsymbol{1}_m \mathbbm{1}_{\Lambda_{n}}(\vect{y},\vect{X})\|
    \sup_{\boldsymbol{\beta}\in\mathcal{B}(B)}\|\boldsymbol{\beta}-\boldsymbol{\beta}_0\| 
    \nonumber\\
    &\leq&  \E\|\mathbb{P}_{n}\varepsilon\vect{X}\boldsymbol{1}_m\mathbbm{1}_{\Lambda_{n}}(\vect{y},\vect{X})\|B.
\end{eqnarray}
And, by Jensen's inequality, we have that
\begin{align*}
    [(\E\|\mathbb{P}_{n}\varepsilon\vect{X}\boldsymbol{1}_m\mathbbm{1}_{\Lambda_{n}}(\vect{y},\vect{X})\|)^2]^{1/2}
    \leq [\E(\|\mathbb{P}_{n}\varepsilon\vect{X}\boldsymbol{1}_m\mathbbm{1}_{\Lambda_{n}}(\vect{y},\vect{X})\|)^2]^{1/2} 
    = [\E\|\vect{X}\boldsymbol{1}_m\|^2\mathbbm{1}_{\Lambda_{n}}(\vect{y},\vect{X})/n]^{1/2}.
\end{align*}
But, from Condition A.2, we have that 
\begin{align*}
    \E\|\vect{X}\boldsymbol{1}_m\|^2\mathbbm{1}_{\Lambda_{n}}(\vect{y},\vect{X}) \leq \E\|\vect{X}\boldsymbol{1}_m\|^2   \leq mC,
\end{align*}
and hence, we finally get
\begin{align*}
    \E(\mathbb{G}_1(B)) \leq 4k_nB(mC/n)^{1/2}.
\end{align*}

Now, on the set $\Lambda_{n}$, we have from Condition A.5 and the Cauchy-Schwarz inequality that
\begin{align*}
    \lvert \rho(\boldsymbol{y}, \boldsymbol{X}\boldsymbol{\beta},{\boldsymbol{\eta}}_0)-\rho(\vect{y},\vect{X}\boldsymbol{\beta}_0,{\boldsymbol{\eta}}_0) \rvert 
    \leq k_n \lvert \boldsymbol{1}_m^T\vect{X}^T(\boldsymbol{\beta}-\boldsymbol{\beta}_0) \rvert  
    \leq k_n\|\vect{X}\boldsymbol{1}_m\| \|\boldsymbol{\beta}-\boldsymbol{\beta}_0\| 
    \leq k_n m^{1/2} K_n B.
\end{align*}
So, we can apply Lemma 4 of \citet{sis_in_generalized_linear_models} with $L^2=4k_n^2mK_n^2B^2$ to get
\begin{align}
\text{Pr}\left(\mathbb{G}_1(B)\geq 4k_nB(mC/n)^{1/2}(1+t)\right) 
    = \exp\big(-2Ct^2/K_n^2\big).
    \label{final0}
\end{align}

\noindent\textit{Step 2:}\\
We generally have that 
\begin{align}
    \text{Pr}(\|\widehat{\boldsymbol{\beta}}-\boldsymbol{\beta}_0\| \geq x) \leq \text{Pr}(\|\widehat{\boldsymbol{\beta}}-\boldsymbol{\beta}_0\| \geq x,  \Omega_n) + \text{Pr}(\Omega_n^c).
\label{EQ:Omegapartition}
\end{align}

Now, in order to find a bound for the first of these two terms, we will in the following concentrate on the set $\Omega_n$ only. Define a convex combination $\boldsymbol{\beta}_s=s\widehat{\boldsymbol{\beta}}+(1-s)\boldsymbol{\beta}_0$, where $ s = (1+\|\widehat{\boldsymbol{\beta}}-\boldsymbol{\beta}_0\|/B)^{-1}$. By definition,
$\|{\boldsymbol{\beta}}_s-\boldsymbol{\beta}_0\| = s \|\widehat{\boldsymbol{\beta}}-\boldsymbol{\beta}_0\|\leq B$,
i.e. ${\boldsymbol{\beta}}_s \in \mathcal{B}(B)$. 

On the set $\Omega_n$ and due to the convexity of $\rho(\cdot)$ in $\boldsymbol{\beta}$ for $\boldsymbol{\beta} \in \boldsymbol{\mathcal{B}}$, we have
\begin{align*}
    \mathbbm{P}_n\rho(\boldsymbol{y}, \boldsymbol{X}\boldsymbol{\beta}_s,\boldsymbol{\eta}_0)-\mathbbm{P}_n\rho(\boldsymbol{y}, \boldsymbol{X}\boldsymbol{\beta}_0,\boldsymbol{\eta}_0) 
    & \leq s[\mathbbm{P}_n\rho(\boldsymbol{y}, \boldsymbol{X}\widehat{\boldsymbol{\beta}},\boldsymbol{\eta}_0)-\mathbbm{P}_n\rho(\boldsymbol{y}, \boldsymbol{X}\boldsymbol{\beta}_0,\boldsymbol{\eta}_0)]\leq 0,
\end{align*}
as $s \in [0,1]$. Next, since $(\boldsymbol{\beta}_0,\boldsymbol{\eta}_0)$ minimizes $\E\rho(\boldsymbol{y},\boldsymbol{X}\boldsymbol{\beta},\boldsymbol{\eta})$, we have 
\begin{align*}
    \E[\rho(\boldsymbol{y}, \boldsymbol{X}\boldsymbol{\beta}_s,\boldsymbol{\eta}_0)-\rho(\boldsymbol{y}, \boldsymbol{X}\boldsymbol{\beta}_0,\boldsymbol{\eta}_0)]\geq 0,
\end{align*}
if we regard ${\boldsymbol{\beta}}_s$ as a parameter in the expectation. Combining the above two results we have that
\begin{align}
    \E[\rho(\boldsymbol{y}, \boldsymbol{X}\boldsymbol{\beta}_s,\boldsymbol{\eta}_0)-\rho(\boldsymbol{y}, \boldsymbol{X}\boldsymbol{\beta}_0,\boldsymbol{\eta}_0)] \leq (\E-\mathbbm{P}_{n})[\rho(\boldsymbol{y}, \boldsymbol{X}\boldsymbol{\beta}_s,\boldsymbol{\eta}_0)-\rho(\boldsymbol{y}, \boldsymbol{X}\boldsymbol{\beta}_0,\boldsymbol{\eta}_0)]
    \leq  \mathbb{G}(B),
\label{EQ:upperboundlogRho}
\end{align}
 where
\begin{align*}
    \mathbb{G}(B)= \sup_{\boldsymbol{\beta}\in \mathcal{B}(B)}\lvert (\mathbb{P}_{n}-E)\{\rho(\boldsymbol{y}, \boldsymbol{X}\boldsymbol{\beta},\boldsymbol{\eta}_0)-\rho(\boldsymbol{y}, \boldsymbol{X}\boldsymbol{\beta}_0,\boldsymbol{\eta}_0)\}\rvert
\end{align*}
By Condition A.4, we have
\begin{align*}
    \|{\boldsymbol{\beta}}_s-\boldsymbol{\beta}_0\| \leq [\mathbb{G}(B)/V]^{1/2}.
\end{align*}
Next, for any $x$, we have that
$\text{Pr}(  \|{\boldsymbol{\beta}}_s-\boldsymbol{\beta}_0\|\geq x)\leq \text{Pr}(  \mathbb{G}(B)\geq V x^2)$.
In particular, letting $x = B/2$, we get 
\begin{align*}
     \text{Pr}(  \|{\boldsymbol{\beta}}_s-\boldsymbol{\beta}_0\|\geq B/2)\leq \text{Pr}(\mathbb{G}(B) \geq VB^2/4).
\end{align*}
By the definition of $\boldsymbol{\beta}_s$, 
$\text{Pr}(\|{\boldsymbol{\beta}}_s-\boldsymbol{\beta}_0\|\geq B/2) = \text{Pr}(\|\widehat{\boldsymbol{\beta}}-\boldsymbol{\beta}_0\|\geq B)$. 
By setting $B=4a_n(1+t)/V$ with $a_n = 4k_n\sqrt{mC/n}$, we get
\begin{align*}
     \text{Pr}(\|\widehat{\boldsymbol{\beta}}-\boldsymbol{\beta}_0\|\geq B) \leq \text{Pr}[\mathbb{G}(B) \geq Ba_n(1+t)],
\end{align*}
and 
\begin{align}
    \text{Pr}[\mathbb{G}(B) \geq Ba_n(1+t)] \leq \text{pr}[\mathbb{G}(B) \geq Ba_n(1+t),\Lambda_n]+\text{Pr}(\Lambda_n^c)
    \label{final1}
\end{align}

Finally, on the set $\Lambda_n$, we have
\begin{align*}
    \sup_{\boldsymbol{\beta}\in \mathcal{B}(B)}\mathbb{P}_{n}\lvert \rho(\boldsymbol{y}, \boldsymbol{X}\boldsymbol{\beta},\boldsymbol{\eta}_0)
    -\rho(\boldsymbol{y}, \boldsymbol{X}\boldsymbol{\beta}_0,\boldsymbol{\eta}_0)\rvert(1-\mathbbm{1}_{\Lambda_{n}}(\vect{y},\vect{X}))=0.
\end{align*}
Hence, by the triangle inequality,
\begin{align*}
\mathbb{G}(B)\leq \mathbb{G}_1(B)+\sup_{\boldsymbol{\beta}\in \mathcal{B}(B)}\lvert\E[ \rho(\boldsymbol{y}, \boldsymbol{X}\boldsymbol{\beta},\boldsymbol{\eta}_0)
-\rho(\boldsymbol{y}, \boldsymbol{X}\boldsymbol{\beta}_0,\boldsymbol{\eta}_0)](1-\mathbbm{1}_{\Lambda_{n}}(\vect{y},\vect{X}))\rvert.
\end{align*}

Then, since Condition A.6 holds for all columns in $\vect{X}$, it follows that 
\eqref{final1} is bounded above by
\begin{align*}
    \text{pr}[\mathbb{G}_1(B) \geq Ba_n(1+t)+o(mC/n)]+ n\text{Pr}(\Lambda_n^c).
\end{align*}

Combining the above with the bound \eqref{final0}, we get  that
\begin{align*}
    \text{Pr}\left(\sqrt{\frac{n}{m}}\left\|\widehat{\boldsymbol{\beta}}-\boldsymbol{\beta}_0\right\|\geq \frac{16k_nC^{1/2}}{V}(1+t),\Omega\right) 
    \leq \exp\left(-\frac{2Ct^2}{K_n^2}\right)+ n\text{Pr}\left(\Lambda_n^c\right).
\end{align*}
Combining the above with \eqref{EQ:Omegapartition}, we get the desired result. \hfill$\square$

\bigskip
\noindent
Note that, 
Theorem \ref{THM:1} gives a bound for the tail probability of the deviation of the M-estimator $\widehat{\boldsymbol{\beta}}$ of the regression coefficients $\boldsymbol{\beta}$ from its true (best-fitting) value $\boldsymbol{\beta}_0$. As a consequence, it also yields the exponential consistency of the M-estimator $\widehat{\boldsymbol{\beta}}$ of the regression coefficient $\boldsymbol{\beta}$
if we additionally assume that the tail probabilities of the response and the covariates go to zero at an exponential rate. We may formally specify these assumptions as follows:

\begin{itemize}
\item[\textit{A.7}] 
There exist positive constants $s_0, s_1$ such that
\begin{align*}
	&\sum_{j=1}^m\Bigg[\E\Big\{\exp\big[b(\vect{x}_{j}^T\boldsymbol{\beta}+\vect{z}_{j}^T\vect{u}+s_0)-b(\vect{x}_{j}^T\boldsymbol{\beta}+\vect{z}_{j}^T\vect{u})\big]\Big\}
	\\
	&~~~~~~~~~~~~~~~+ \E\Big\{\exp\big[b(\vect{x}_{j}^T\boldsymbol{\beta}+\vect{z}_{j}^T\vect{u}-s_0)-b(\vect{x}_{j}^T\boldsymbol{\beta}+\vect{z}_{j}^T\vect{u})\big]\Big\}\Bigg] \leq s_1.
\end{align*}
 \item[\textit{A.8}]
    There exist positive constants $r_0, r_1, \delta$ such that for sufficiently large $t$, we have
	\begin{align*}
	\text{Pr}(\| \vect{X}\|_\infty > t) \leq (r_1-s_1)\exp(-r_0t^\delta),
\end{align*}
where $s_1$ is the same constant as in Condition A.7. 
 \item[\textit{A.9}]
There exist constants $r_2$ and $r_3$ such that
$\text{Pr}(\Omega_n^c) \leq r_2\exp(-r_3n)$
\end{itemize}

Assumption A.7 requires that the function $b(\cdot)$ in the exponential family in \eqref{EQ:exponentialfamily} does not increase or decrease too fast. A.8 assumes a strong tail decay of the components of $\vect{X}$, which holds for e.g. all sub-Gaussian distributions with $r_1-s_1=2$ and $\delta = 2$, though A.8 is more general. 
Assumptions A.7--A.8 are generalizations of Condition D in \cite{sis_in_generalized_linear_models}; these types of assumptions are commonly used in the literature on variable screening 
(see e.g. \cite{abhikponzi} and \cite{conditional-screening}) and ensure that the response variable $\vect{y}$ has an exponentially light tail, as shown in the proof of Theorem 2 below. 
Assumption A.9 is somewhat similar to A.3, in the sense that the empirical mean of the loss function should be convex (in both $\boldsymbol{\beta}$ and $\boldsymbol{\eta}$) in a neighborhood around $(\boldsymbol{\widehat{\beta}},\boldsymbol{\widehat{\eta}})$, asymptotically. 
It can be verified that Assumption A.9 holds if  $\boldsymbol{\widehat{\eta}}$ is consistently estimated at an appropriate rate and $\rho$ is well-behaved around $(\boldsymbol{\beta}_0, \boldsymbol{\eta}_0)$.
Under these additional assumptions, we have the exponential consistency of M-estimators of regression coefficients under GLMMs as stated in the following theorem. 

\begin{theorem}[Exponential consistency]
Suppose that Assumptions A.1 -- A.9 hold under the set-up of generalized linear mixed models and a given $\rho$-function defining the M-estimator of the parameters as in \eqref{EQ:M-est}. 
Then, as $n\rightarrow\infty$ with $n/(K_n^2k_n^2m)\rightarrow \infty$, where $K_n$ and $k_n$ are defined as in A.5, 
the M-estimator $\widehat{\boldsymbol{\beta}}$ of the regression coefficient $\boldsymbol{\beta}$ converges to its true (best-fitting) value $\boldsymbol{\beta}_0$, in probability,
at an exponential rate.
\end{theorem}
\textit{Proof:}
By the exponential Chebyshev's inequality and the law of total expectation, we have that 
\begin{align*}
    \text{Pr}(\lvert y_j \lvert \geq u)  = \text{Pr}(y_j \geq u) + \text{Pr}(-y_j \geq u)\leq \exp(-s_0u)\E[\exp(s_0y_j)+\exp(-s_0y_j)] \\=  \exp(-s_0u)\E\{\E[\exp(s_0y_j)\lvert \theta_j]+\E[\exp(-s_0y_j)\lvert \theta_j]\}.
\end{align*}
Since $y_j$ given $\theta_j$ belongs to the exponential family, we have that
\begin{align*}
    \E[\exp(s_0y_j)\lvert \theta_j]\ &= \exp\{b(\theta_j+s_0)-b(\theta_j)\}, \quad \text{and} \\
    \E[\exp(-s_0y_j)\lvert \theta_j]\ &= \exp\{b(\theta_j-s_0)-b(\theta_j)\}
\end{align*}

By Condition A.7, we then get
\begin{align}
    \text{Pr}(\|\vect{y}\|_\infty \geq u) &\leq \sum_{j=1}^m   \text{Pr}(\lvert y_j\rvert \geq u) \nonumber \\
    &\leq \exp(-s_0u)\sum_{j=1}^m\Bigg[\E\Big\{\exp\big[b(\vect{x}_{j}^T\boldsymbol{\beta}+\vect{z}_{j}^T\vect{u}+s_0)-b(\vect{x}_{j}^T\boldsymbol{\beta}+\vect{z}_{j}^T\vect{u})\big]\Big\}
	\nonumber \\&+\E\Big\{\exp\big[b(\vect{x}_{j}^T\boldsymbol{\beta}+\vect{z}_{j}^T\vect{u}-s_0)-b(\vect{x}_{j}^T\boldsymbol{\beta}+\vect{z}_{j}^T\vect{u})\big]\Big\}\Bigg]\nonumber \\
& \leq s_1\exp(-s_0u).
 \label{EQ:upperbound||y||}
\end{align}
Now, we use this result to show that $\text{Pr}(\Lambda_n^c)$ in Theorem 1 vanishes at an exponential rate. By Condition A.8 and by setting $u = r_0K_n^{\delta}/s_0$ in \eqref{EQ:upperbound||y||}, we have that
\begin{align*}
    \text{Pr}(\Lambda_n^c) \leq \text{Pr}(\|\vect{X}\|_\infty > K_n) + \text{Pr}(\|\vect{y}\|_\infty > r_0K_n^{\delta}/s_0) \leq r_1\exp(-r_0K_n^\delta),
\end{align*}
where $K_n \to \infty$ as $n \to \infty$. By taking $1+t = c_1V_nn^{1/2}/(16k_nC^{1/2}m^{1/2})$ for some constant $c_1$ and using Assumption A.9, we get that 
\begin{align*}
\text{Pr}\left(\left\|\widehat{\boldsymbol{\beta}}-\boldsymbol{\beta}_0\right\|\geq c_1\right) 
    \leq \exp\left(-\frac{c_2}{K_n^2k_n^2}\frac{n}{m}\right)+ nr_1\exp(-r_0K_n^\delta)+r_2\exp(-r_3n),
\end{align*}
for some constant $c_2$ and $n/(K_n^2k_n^2m) \to \infty$ as $n \to \infty$. This completes the proof of Theorem 2.  $\hfill \square$

\section{Illustrations}

In this section, we will revisit the examples of M-estimators introduced in Section \ref{sec:Mestforglmm} in order to justify the conditions required for Theorem 1 and 2 under two specific cases of generalized linear mixed regression settings, namely the linear mixed model and the logistic mixed model, with theoretical and/or empirical illustrations of their exponential consistency. In this respect, we note that Condition A.1, A.2 and A.8 are generic ones, with A.2 and A.8  being dependent on the boundedness of the covariates, and hence they will always assumed to be true and will not be discussed for the following specific cases to avoid repetitions.

\subsection{Linear Mixed Effects Regression}
\label{Linearmixedeffects}
Consider the  linear mixed effects model, the simplest case of a GLMM, defined as  
\begin{align*}
\vect{y} = \vect{X}\boldsymbol{\beta}+\vect{Z}\vect{u}+\boldsymbol{\epsilon},
\end{align*}
where $\boldsymbol{\epsilon} \sim \mathcal{N}(0, \sigma_0^2\mathbb{I}_m)$ and $\vect{X}$, $\vect{Z}$ and $\vect{u}$ are defined as in Section \ref{sec:Mestforglmm}. Assume we have observations $(\vect{y}_i,\vect{X}_i)_{i=1,...,n}$ as independent realizations of $(\vect{y},\vect{X})$. The link function in \eqref{EQ:conditionalmeanofy} is the identity, i.e.  $h^{-1}(\theta_j) = \theta_j$, so the marginal expected value of $\vect{y}_i$ is $\E(\vect{y}_i) =\E[\E(\vect{y}_i\lvert \vect{u})]  = \vect{X}_i\boldsymbol{\beta}$, for $i = 1,...,n$, where we again have omitted the the explicit conditioning on $\vect{X}$ for simplicity. We define the variance parameter vector $\boldsymbol{\eta} = (\sigma_0^2,\sigma_1^2,...,\sigma_{q(q+1)/2}^2)$. In this setting, we know the marginal distribution of each response vector $\vect{y}_i$:
\begin{align}
    \vect{y}_i \sim \mathcal{N}(\vect{X}_i\boldsymbol{\beta},\vect{V}_i) \quad i = 1,...,n, \quad \text{ where } \vect{V}_i = \vect{V}_i(\boldsymbol{\eta})=\sigma_0^2\mathbb{I}_{m}+\vect{Z}_{i}\vect{G}\vect{Z}_{i}^T.
 \label{EQ:yNormallydistributed}
\end{align}
We now make the following simplifying assumptions which are needed for Assumptions A.3--A.7 to hold in the linear mixed models setting:
\begin{itemize}
\item[\textit{B.1}] 
The matrix $\vect{V}$ is positive definite and $\vect{X}$ is of full column rank 
\item[\textit{B.2}]
Assume that there exist positive constants $s_0$ and $s_3$ such that 
\begin{align*}
   \E\Big[\exp\{s_0(\vect{x}_j^T\boldsymbol{\beta}+\vect{z}_j^T\vect{u})\}+\exp\{-s_0(\vect{x}_j^T\boldsymbol{\beta}+\vect{z}_j^T\vect{u})\}\Big] \leq s_3, \quad j=1,...,m.
\end{align*}
\end{itemize}

\subsubsection*{Example: MLE}
Since we know the marginal distribution of the response, we have a simple closed form expression for the loss function associated with the MLE, which is the negative log-likelihood function having the form 
\begin{align*}
\rho(\vect{y},\vect{X}\boldsymbol{\beta},\boldsymbol{\eta}) = \frac{m}{2}\log(2\pi)+\frac{1}{2}\log\lvert \vect{V}({\boldsymbol{\eta}})\rvert+\frac{1}{2}(\vect{y}-\vect{X}\boldsymbol{\beta})^T\vect{V}({\boldsymbol{\eta}})^{-1}(\vect{y}-\vect{X}\boldsymbol{\beta}).
\end{align*}
By differentiating this expression with respect to $\boldsymbol{\beta}$, we get the estimating equations for estimating $\boldsymbol{\beta}$ as given by
\begin{align*}
\mathbbm{P}_n\frac{\partial}{\partial \boldsymbol{\beta}}\rho(\vect{y},\vect{X}\boldsymbol{\beta},\boldsymbol{\eta}) = \frac{1}{n}\sum_{i}^n
\vect{X}_i^T\vect{V}_i^{-1}(\vect{y}_i-\vect{X}_i\boldsymbol{\beta}) = 0.
\end{align*}
Similarly, we estimate $\boldsymbol{\eta}$ by solving the estimating equations for each $\sigma_r^2$, $r = 0,...,\frac{q(q+1)}{2}$:
\begin{align}
\mathbbm{P}_n\frac{\partial}{\partial \sigma_r^2}\rho(\vect{y},\vect{X}\boldsymbol{\beta},\boldsymbol{\eta}) = \frac{1}{n}\sum_{i}^n \Big\{(\vect{y}_i-\vect{X}_i\boldsymbol{\beta})^T\vect{V}_i^{-1}\vect{U}_{ir}\vect{V}_i^{-1}(\vect{y}_i-\vect{X}_i\boldsymbol{\beta})-\text{Tr}\big(\vect{V}_i^{-1}\vect{U}_{ir}\big)\Big\}=0,
\label{EQ:MLEestimatingeqsSigma}
\end{align}
where $\vect{U}_{ir}$ is defined as the partial derivative of $\vect{V}_i$ with respect to $\sigma_r^2$. We now argue that the MLE for the linear mixed model satisfies the regularity conditions for Theorem 1 under the simplified intuitive assumptions B.1 and B.2. 
\begin{lemma}
  Suppose that Assumptions B.1 and B.2 hold under the set-up of linear mixed models and the function $\rho(\cdot)$ is defined as the negative log-likelihood. Then, Assumptions A.3--A.7 hold. 
\end{lemma}

The proof of Lemma 3 is given in Appendix A. Now, if we assume that A.1--A.2 and A.8--A.9 hold in addition to the assumptions in Lemma 3, Theorem 1 and 2 will follow directly for the MLE for linear mixed models.

\subsubsection*{Example: MDPDE}

The MDPDE was formulated in the setting of linear mixed models by \cite{saraceno2020robust} and the loss function for the MDPDE is given by (which can also be obtained by simplifying \eqref{EQ:MDPDEobjfunc} directly for this case)
\begin{align}
\rho(\vect{y},\vect{X}\boldsymbol{\beta},\boldsymbol{\eta}) = \frac{1}{(2\pi)^{m\alpha/2}\lvert\vect{V}\rvert^{\alpha/2}(1+\alpha)^{m/2}}-\Big(1+\frac{1}{\alpha}\Big)\frac{1}{(2\pi)^{m\alpha/2}\lvert\vect{V}\rvert^{\alpha/2}}\exp\big\{-\frac{\alpha}{2}(\vect{y}-\vect{X}\boldsymbol{\beta})^T\vect{V}^{-1}(\vect{y}-\vect{X}\boldsymbol{\beta})\big\}\\ 
    = L_1\lvert\vect{V}\rvert^{-\alpha/2}-L_2\lvert\vect{V}\rvert^{-\alpha/2}\exp\Big\{-\frac{\alpha}{2}(\vect{y}-\vect{X}\boldsymbol{\beta})^T\vect{V}^{-1}(\vect{y}-\vect{X}\boldsymbol{\beta})\Big\},
    \label{EQ:MDPDElossfunctionLinearmixed}
\end{align}
where $L_1 = (1+\alpha)^{-m/2}(2\pi)^{-m\alpha/2}$ and $L_2 = \Big(1+\frac{1}{\alpha}\Big)(2\pi)^{-m\alpha/2}.$     Differentiating with respect to $\boldsymbol{\beta}$, we get the estimating equation for the MDPDE of $\boldsymbol{\beta}$ as
\begin{align}
    \mathbb{P}_n\frac{\partial \rho(\boldsymbol{y}, \boldsymbol{X}\boldsymbol{\beta},\boldsymbol{\eta})}{\partial \boldsymbol{\beta}} = \frac{1}{n}\sum_{i=1}^n\Big[-\Big(1+\frac{1}{\alpha}\Big)\frac{e^{-\frac{\alpha}{2}(\vect{y}_i-\vect{X}_i\boldsymbol{\beta})^T\vect{V}_i^{-1}(\vect{y}_i-\vect{X}_i\boldsymbol{\beta})}}{2\pi^{\frac{n\alpha}{2}}\lvert \vect{V}_i\rvert^{\frac{\alpha}{2}}}\alpha\vect{X}_i^T\vect{V}_i^{-1}(\vect{y}_i-\vect{X}_i\boldsymbol{\beta})\Big] = 0.
    \label{EQ:mdpdeforBetaEstimating}
 \end{align}
The estimating equations with respect to each element in $\boldsymbol{\eta}$ is similarly obtained as
\begin{align}
   \mathbb{P}_n\frac{\partial \rho(\boldsymbol{y}, \boldsymbol{X}\boldsymbol{\beta},\boldsymbol{\eta})}{\partial \sigma_r^2} = \frac{1}{n}\sum_{i=1}^n\Bigg\{-\frac{\alpha\text{Tr}(\vect{V}_i^{-1}\vect{U}_{ir})}{2(2\pi)^{\frac{n\alpha}{2}}\lvert\vect{V}_i \rvert^{\frac{\alpha}{2}}(\alpha+1)^{\frac{m}{2}}}+\frac{\alpha}{2}\Big(1+\frac{1}{\alpha}\Big)\frac{e^{-\frac{\alpha}{2}(\vect{y}_i-\vect{X}_i\boldsymbol{\beta})^T\vect{V}_i^{-1}(\vect{y}_i-\vect{X}_i\boldsymbol{\beta})}}{(2\pi)^{\frac{n\alpha}{2}}\lvert \vect{V}_i\rvert^{\frac{\alpha}{2}}}\nonumber\\ 
    \times \Big[\text{Tr}(\vect{V}_i^{-1}\vect{U}_{ir})-(\vect{y}_i-\vect{X}_i\boldsymbol{\beta})^T\vect{V}_i^{-1}\vect{U}_{ir}\vect{V}_i^{-1}(\vect{y}_i-\vect{X}_i\boldsymbol{\beta})\Big]\Bigg\} = 0,
\label{EQ:mdpdeforEtaEstimating}
\end{align}
where $\vect{U}_{ir}$ is defined as in \eqref{EQ:MLEestimatingeqsSigma}. We may note here that these estimating equations are equivalent to the estimating equations of the MLE by letting $\alpha \to 0$; this is because the MDPDE coincides with the MLE at $\alpha=0$ and provide a robust generalization of the MLE at $\alpha>0$ with increasing robustness with increase in $\alpha>0$ (see \cite{saraceno2020robust} for details). Let $\vect{d} := \vect{V}^{-1}(\vect{y}-\vect{X}\boldsymbol{\beta})$. We then have the following simplifying assumptions that are needed for A.3--A.7 to hold for the MDPDE in linear mixed models:

\begin{itemize}
\item[\textit{B.3}] For any $\vect{v}\in \mathbb{R}^p$ and $\boldsymbol{\beta} \in \boldsymbol{\mathcal{B}}$, $\vect{v}^T\vect{X}^T\vect{V}^{-1}\vect{X}\vect{v} \geq \alpha\vect{v}^T\vect{X}^T\vect{d}\vect{d}^T\vect{X}\vect{v}$.
\item[\textit{B.4}] The matrix $$\vect{V}^{-1}\otimes\vect{V}^{-1}\text{vec}(\vect{V})\text{vec}(\vect{V}^{-1})^T$$
is positive semi-definite, where $\text{vec}(\vect{V})$ is  the $m^2 \times 1$ column vector obtained by stacking the columns of the matrix $\vect{V}$, and $\otimes$ denotes the Kronecker product.
\item[\textit{B.5}] The matrix $\E[e^{-\alpha (\vect{y}-\vect{X}\boldsymbol{\beta})^T\vect{V}^{-1}(\vect{y}-\vect{X}\boldsymbol{\beta})}\vect{d}\vect{d}^T\otimes\vect{d}\vect{d}^T]$ is positive definite at $(\boldsymbol{\beta}_0,\boldsymbol{\eta}_0)$.
\end{itemize}

\begin{lemma}
    Suppose that B.1--B.5 hold under the set-up of linear mixed models and the function $\rho(\cdot)$ is defined as the minimum density power divergence in \eqref{EQ:mdpdelossfunctionGeneral}. Then Assumptions A.3--A.7 hold.
\end{lemma}
Proof is given in appendix B. Hence, given the conditions for Lemma 4, Theorem 1 and 2 follows directly for the MDPDE in the linear mixed model setting under Assumptions A.1--A.2, A.8--A.9 and B.1--B.5.

\subsubsection{Empirical Illustration}
\label{sec:linearmixedsimulation}

\begin{figure}[H]
\centering
\includegraphics[width=0.75\linewidth]{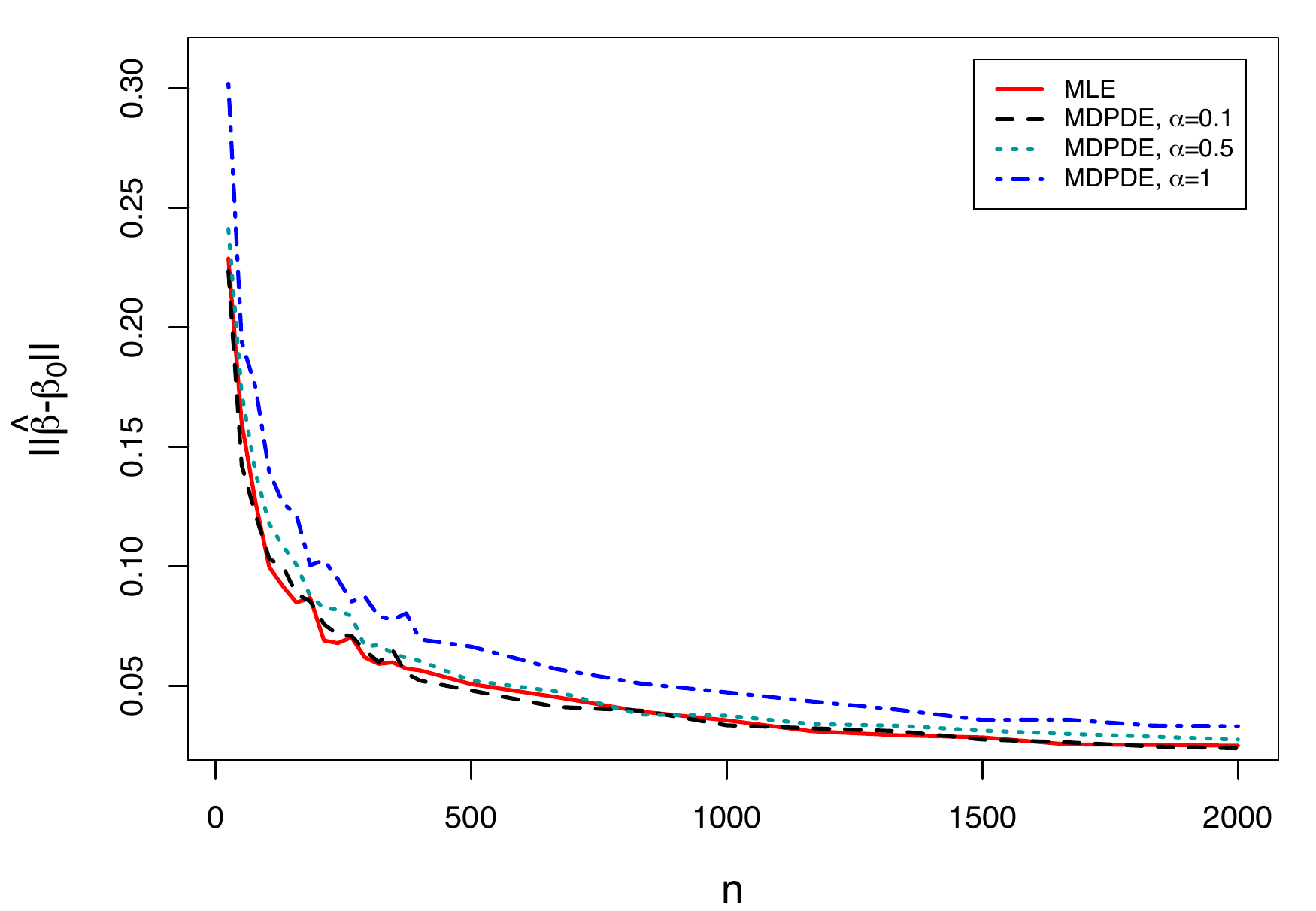}
\caption{The $\ell_2$ distance between the estimate $\widehat{\boldsymbol{\beta}}$ and the population parameter $\boldsymbol{\beta}_0$ as a function of the sample size $n$, for MLE and MDPDE for $\alpha \in \{0.1,0.5,1\}$ in the linear mixed model example, based on 150 replications.}
\label{fig:fig}
\end{figure}

Let us now empirically validate our exponential consistency result under the linear mixed model through simulations, by showing that  $\|\widehat{\boldsymbol{\beta}}-\boldsymbol{\beta}_0\|$ indeed decreases exponentially 
as a function of $n$ for both the MLE and the MDPDE with different $\alpha>0$. We use a similar simulation set-up as in \cite{ghoshthoresen2018}. We let the number of groups (independent observations) $n$ vary, and set the number of observations from each group to $m = 6$. We take $p=5$ with the true parameter being $\boldsymbol{\beta}_0 = (1,2,4,3,3)^T$. The number of random effects is taken to be $q=2$, with the random effects vector $\vect{u} \sim \mathcal{N}(0,\sigma_u^2\mathbb{I}_2)$, i.e. the random effects are independent. We set $\sigma_u^2=0.56$ and the error variance $\sigma^2=0.25$. The first column in the design matrix $\vect{X}_i$ is equal to 1 (the intercept), while the rest of the columns are sampled from a multivariate normal distribution with mean 0 and identity covariance matrix. The design matrix of the random effects,  $\vect{Z}_i$ is the first two columns of $\vect{X}_i$ for each group $i=1,...,n$. Figure \ref{fig:fig} shows the average estimated $\ell_2$ bias $\|\widehat{\boldsymbol{\beta}}-\boldsymbol{\beta}_0\|$ as a function of $n$, based on 150 replications, for the MLE and the MDPDE with different values of $\alpha>0$. It can be clearly seen from the figure that for both MLE and MDPDE, the rate of convergence is indeed exponential, with the rate of convergence being similar for all the estimators considered, although the exact values of their biases differ slightly. Consistent with the literature of the MDPDE, the bias under pure data increases slightly  with increasing values of $\alpha$ (a price to pay for greater robustness), but what is new here is that their rate of convergence are interestingly similar (exponential) for all $\alpha\geq 0$.

\subsection{Logistic mixed effects regression}
We now consider a binary mixed effects model. Similar to Section \ref{Linearmixedeffects}, we denote the fixed and random effects design matrices by $\vect{X}$ and $\vect{Z}$, and the vectors of fixed and random effects coefficients by $\boldsymbol{\beta}$ and $\vect{u}$, respectively. We assume $\vect{u}\sim \mathcal{N}(0,\vect{G}(\boldsymbol{\eta}))$. Consider observations from $n$ independent groups, i.e. we have $n$ independent response vectors $\vect{y}_i = [y_{i1},...,y_{im}]^T$. The logistic mixed effects model for clustered binary data is then specified through the mean function $\mu_{ij}$ for the response vector $\vect{y}_i$:
\begin{align}
   \mu_{ij} = \E[y_{ij}\lvert \vect{u}] = \frac{\exp(\vect{x}_{ij}^T\boldsymbol{\beta}+\vect{z}_{ij}^T\vect{u})}{1+\exp(\vect{x}_{ij}^T\boldsymbol{\beta}+\vect{z}_{ij}^T\vect{u})} \quad \text{ and } \quad
    \Var[y_{ij}\lvert \vect{u}] = \frac{\exp(\vect{x}_{ij}^T\boldsymbol{\beta}+\vect{z}_{ij}^T\vect{u})}{[1+\exp(\vect{x}_{ij}^T\boldsymbol{\beta}+\vect{z}_{ij}^T\vect{u})]^2},
    \label{EQ:defineLogisticMixed}
\end{align}
with the conditional distribution of $y_{ij}$ being
\begin{align}
    f_{y_{ij}|\vect{u}}(y_{ij}\lvert \vect{u},\boldsymbol{\beta}) = \mu_{ij}^{y_{ij}}(1-\mu_{ij})^{1-y_{ij}}, \quad  j = 1,..., m.
    \label{EQ:condYlogistic}
\end{align}

\subsubsection*{Example: MLE}
Consider now the random response vector $\vect{y} =[y_1,...,y_m]^T$, of which $(\vect{y}_1,...,\vect{y}_n)$ are independent realizations. The negative log-likelihood function for the logistic mixed effects model takes the form
\begin{align*}
\rho(\vect{y},\vect{X}\boldsymbol{\beta},\boldsymbol{\eta}) = \frac{mq}{2}\log(2\pi) + \frac{m}{2}\log\lvert\vect{G}(\boldsymbol{\eta})\rvert+\sum_{j=1}^m\log\int \exp\{h_j(y_j,\vect{u};\boldsymbol{\beta},\boldsymbol{\eta})\}d\vect{u},
\end{align*}
where
\begin{align}
h_j(y_j,\vect{u};\boldsymbol{\beta},\boldsymbol{\eta})= y_j(\vect{x}_{j}^T\boldsymbol{\beta}+\vect{z}_j^T\vect{u})-0.5\vect{u}\vect{G}(\boldsymbol{\eta})^{-1}\vect{u}-\log\big[1+\exp(\vect{x}_{j}^T\boldsymbol{\beta}+\vect{z}_{j}^T\vect{u})\big],
\label{EQ:logisticgaussianintegral}
\end{align}

As before, the MLE for $\boldsymbol{\beta}$ and $\boldsymbol{\eta}$ is to be obtained as the solution to the score equations $\mathbbm{P}_n\frac{\partial}{\partial\boldsymbol{\beta}}\rho(\vect{y},\vect{X}\boldsymbol{\beta},\boldsymbol{\eta}) = 0$ and $\mathbbm{P}_n\frac{\partial}{\partial\boldsymbol{\eta}}\rho(\vect{y},\vect{X}\boldsymbol{\beta},\boldsymbol{\eta}) = 0$, which need to be solved numerically, with some approximation of the integral. This also means that we need some approximations to compute the information matrix. Equation (7.62) in \cite{demidenko2013mixed} gives an approximation to the Fisher information matrix for $\boldsymbol{\beta}$ in the case of a single random effect and shows that this approximation is positive definite. In the case of multiple random effects, we get a multidimensional integral and the approximation becomes more complicated, in particular when we want to find the full information matrix. Thus, we will study the exponential consistency of the MLE under the logistic mixed-models through simulation studies in Section \ref{sec:logisticsimulation}.

\subsubsection*{Example: MDPDE}
For the logistic mixed effects model, the MDPDE can be defined by the minimizer of the general loss function described in \eqref{EQ:MDPDEobjfunc} but with the specific form of the density $f(\vect{y};\boldsymbol{\beta},\boldsymbol{\eta})$ as given by    
\begin{align*}
f(\vect{y};\boldsymbol{\beta},\boldsymbol{\eta}) = \prod_{j=1}^m\int \mu_{j}^{y_{j}}(1-\mu_{j})^{1-y_{j}}f_{\vect{u}}(\vect{u}\lvert \boldsymbol{\eta})d\vect{u} 
\end{align*}
and $f_{\vect{u}}(\vect{u}\lvert \boldsymbol{\eta})$ as in \eqref{EQ:mixedeffectsdistribution}. Therefore, the MDPDE loss function for the present case can be simplified as 

\begin{align*}
\rho(\vect{y},\vect{X}\boldsymbol{\beta},\boldsymbol{\eta})= \frac{1}{(2\pi)^{(1+\alpha)mq/2}\lvert\vect{G}\rvert^{(1+\alpha)m/2}}\sum_{\vect{y}\in \{0,1\}^m}\prod_{j=1}^m\Big[\int \exp\{h_j(y_j,\vect{u};\boldsymbol{\beta},\boldsymbol{\eta})\}d\vect{u}]^{1+\alpha} \\-\big(1+\frac{1}{\alpha}\big)\prod_{i=1}^m\Big[\int \exp\{h_j(y_j,\vect{u};\boldsymbol{\beta},\boldsymbol{\eta})\}d\vect{u}]^{\alpha},
\end{align*}
where $h_j(y_j,\vect{u};\boldsymbol{\beta},\boldsymbol{\eta})$ is given in \eqref{EQ:logisticgaussianintegral}. For the same reasons as the MLE, the information matrix and the Hessian of the loss function is not possible to analyze analytically, so we will explore its properties through simulations below.

\subsubsection{Empirical Illustration}
\label{sec:logisticsimulation}

\begin{figure}[H]
\centering
\includegraphics[width=0.75\linewidth]{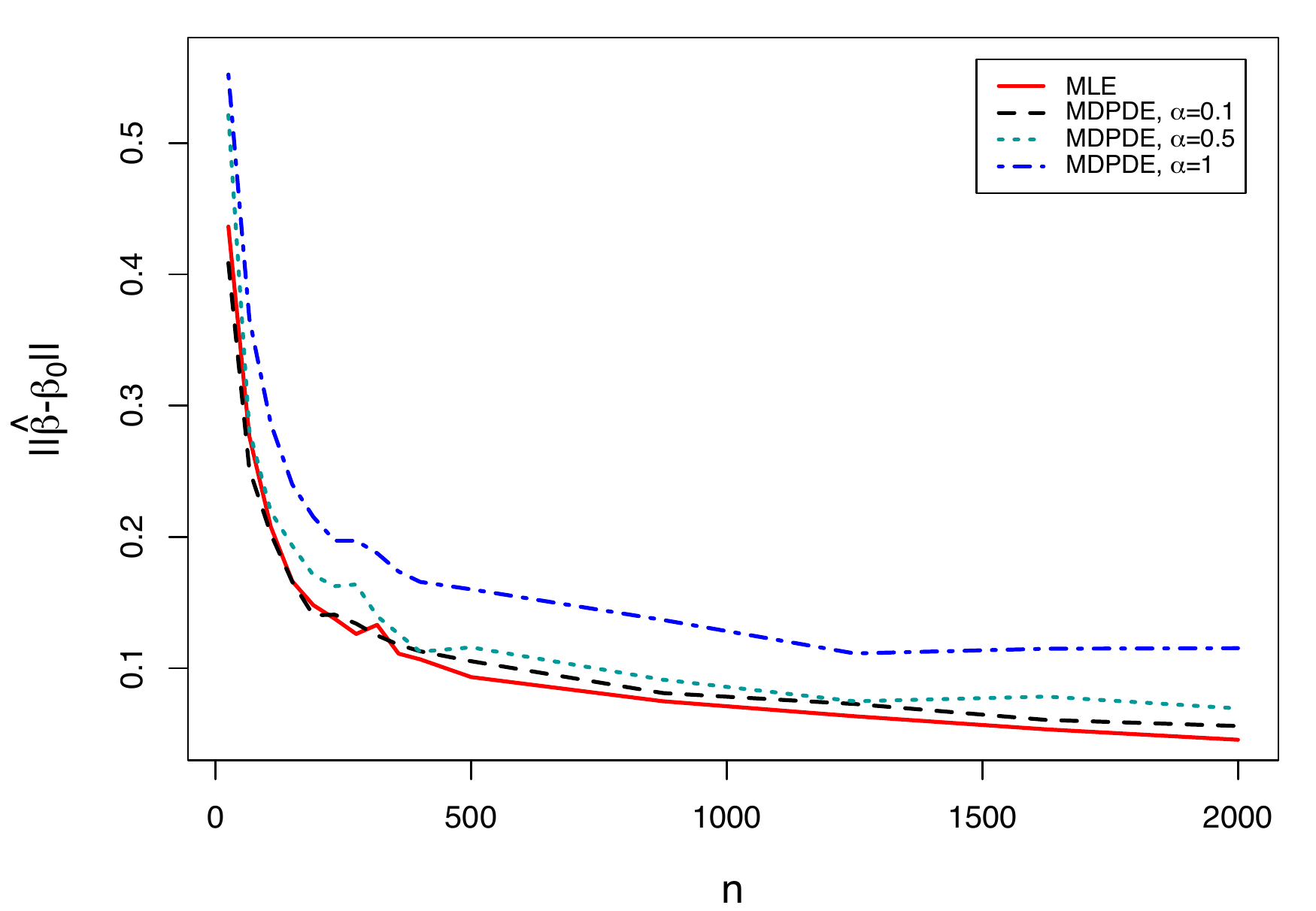}
\caption{The $\ell_2$ distance between the estimate $\widehat{\boldsymbol{\beta}}$ and the population parameter $\boldsymbol{\beta}_0$ as a function of the sample size $n$, for MLE and MDPDE for $\alpha \in \{0.1,0.5,1\}$ in the logistic mixed model example, based on 150 replications.}
\label{fig:fig2}
\end{figure}

As it is difficult to verify the theoretical assumptions directly for the present set-up, we will empirically  justify that our exponential consistency results indeed also hold for logistic mixed models, through simulations. In particular, let us consider the simulation set-up similar to the LMM example described in Section \ref{sec:linearmixedsimulation}, with the number of observations within each group being $m = 6$, but now we set the number of fixed effects $p=2$, and $\boldsymbol{\beta}_0 = (1,2)^T$. The number of random effects is $q=1$, with $u \sim \mathcal{N}(0,\sigma_u^2)$, and $\sigma_u^2=0.56$. The design matrices $\vect{X}$ and $\vect{Z}$  are also simulated similar to the ones described in Section 4.1.1., i.e. we have a random intercept model.  The response for group $i$ and measurement $j$ is then generated following $y_{ij} \lvert \vect{u} \sim \text{Binom}(\mu_{ij})$, where $\mu_{ij}$ is given by \eqref{EQ:defineLogisticMixed}. Figure \ref{fig:fig2} shows the plot of $\|\widehat{\boldsymbol{\beta}}-\boldsymbol{\beta}_0\|$ as a function of $n$, averaged over 150 replications, for the MLE and the MDPDE with different $\alpha>0$. It can again be seen from the figure that the rate of convergence of this $\ell_2$ bias to zero is exponential, reconciling the fact that our theoretical exponential consistency results continue to hold for such logistic mixed models as well.

\section{Conclusions}

In this paper, we have proved that, for the class of GLMMs, M-estimators with sufficiently smooth loss functions will yield regression coefficient estimates that converge to the true value in probability, at an exponential rate. Importantly, these results show that the convergence in probability is, in general, faster for M-estimators than previously known results under the GLMM set-up. In addition, our results hold without assuming that the variance parameter $\boldsymbol{\eta}$ is known.
Through simulation studies, we learnt that the rate of convergence is similar for the classical (non-robust) MLE and the robust MDPDE for different $\alpha>0$, which is also an addition to the literature of the popular MDPDE as a robust generalization of the MLE.  

Similar to \cite{sis_in_generalized_linear_models}, our exponential consistency results are also of importance in the setting of high-dimensional regression. Theorem 1 and 2, together with some additional conditions, may be utilized to prove the sure screening property of variable screening procedures based on M-estimates of marginal regression coefficients, indicating that the probability of capturing all truly important covariates will tend to one for such procedures besides some of them being robust against data contamination as well. This makes it possible to define a wide range of variable screening procedures under ultra-high dimensional GLMMs that are both robust and have the sure screening property. We hope to explore this avenue further to study such robust variable screening procedures for GLMMs in our future research. 

Furthermore, extending our exponential consistency results for other types of existing robust estimators under the GLMM set-up (e.g., those proposed by \cite{sinha2004robust} and \cite{copt2006high}) and also for M-estimators under more general multivariate models, possibly with a matrix of regression coefficients instead of vectors, would be important and we hope to consider this in our future research work as well.

\bigskip\bigskip

\bibliography{references}

\newpage

\appendix

\section{Proof of Lemma 3}

For Condition A.3, the information matrix of the loss function is given by Eq. (6.62) in \cite{searlemcculloch}:
\begin{align*}
    \boldsymbol{I}(\boldsymbol{\beta}, \boldsymbol{\eta}) 
    = \E\begin{bmatrix}\vect{X}^T\vect{V}^{-1}\vect{X} & \vect{0}\\ 
    \vect{0} & \frac{1}{2}[\text{tr}(\vect{Z}_j^T\vect{V}^{-1}\vect{Z}_k(\vect{Z}_j^T\vect{V}^{-1}\vect{Z}_k)^T)]_{0\leq j, k\leq 1}\end{bmatrix},
\end{align*}
where $\vect{Z}_0 := \mathbb{I}_m$. This information matrix is finite and positive definite for all $\boldsymbol{\beta}$ and $\boldsymbol{\eta}$ by Assumption B.1.

 This also implies that $\|\boldsymbol{I}(\boldsymbol{\beta}, \boldsymbol{\eta}) \|$ is bounded from above by the maximum eigenvalue of the information matrix. The expression $\vect{X}^T\vect{V}^{-1}\vect{X}$ is the Hessian of the loss function with respect to $\boldsymbol{\beta}$ and is positive definite by the same arguments. Thus, the loss function is convex in $\boldsymbol{\beta}$. Hence, Condition A.3 holds.
By the convexity of $\rho(\cdot)$ in $\boldsymbol{\beta}$, condition A.4 holds because
\begin{align*}
    \E[\rho(\vect{y},\vect{X}\boldsymbol{\beta},\boldsymbol{\eta}_0)-\rho(\vect{y},\vect{X}\boldsymbol{\beta}_0,\boldsymbol{\eta}_0)] \geq \E\big[\nabla_\beta\rho(\vect{y},\vect{X}\boldsymbol{\beta}_0,\boldsymbol{\eta}_0)^T(\boldsymbol{\beta}-\boldsymbol{\beta}_0)\big]\\
    = \E\big[(\vect{X}^T\vect V(\boldsymbol{\eta}_0)^{-1}\vect{y}-\vect{X}^T\vect V(\boldsymbol{\eta}_0)^{-1}\vect{X}\boldsymbol{\beta}_0)^T(\boldsymbol{\beta}-\boldsymbol{\beta}_0)\big]\\
    = [\vect{X}^T\vect V(\boldsymbol{\eta}_0)^{-1}\vect{X}\boldsymbol{\beta}-\vect{X}^T\vect V(\boldsymbol{\eta}_0)^{-1}\vect{X}\boldsymbol{\beta}_0)]^T(\boldsymbol{\beta}-\boldsymbol{\beta}_0) \\= (\boldsymbol{\beta}-\boldsymbol{\beta}_0)^T\vect{X}^T\vect V(\boldsymbol{\eta}_0)^{-T}\vect{X}(\boldsymbol{\beta}-\boldsymbol{\beta}_0) \geq V\|\boldsymbol{\beta}-\boldsymbol{\beta}_0\|^2
\end{align*}
where $V$ is the minimum eigenvalue of the matrix $\vect{X}^T\vect V(\boldsymbol{\eta}_0)^{-T}\vect{X}$, which is positive since the matrix is positive definite. Condition A.5 holds as the function $\rho(\vect{y},\vect{X}\boldsymbol{\beta},\boldsymbol{\eta})$ is continuously differentiable and locally bounded in $\boldsymbol{\mathcal{B}}$, and thus its gradient is locally bounded. 
For condition A.6, we have that
\begin{align}
 	\E[\rho(\boldsymbol{y}, \boldsymbol{X}\boldsymbol{\beta},\boldsymbol{\eta}_0)-\rho(\boldsymbol{y}, \boldsymbol{X}\boldsymbol{\beta}_0,\boldsymbol{\eta}_0)]
	(1-\mathbbm{1}_{\Lambda_{n}}(\vect{y},\vect{X})) =  \E[\rho(\boldsymbol{y}, \boldsymbol{X}\boldsymbol{\beta},\boldsymbol{\eta}_0)-\rho(\boldsymbol{y}, \boldsymbol{X}\boldsymbol{\beta}_0,\boldsymbol{\eta}_0)] \nonumber \\ 
	- \{\E[ \rho(\boldsymbol{y}, \boldsymbol{X}\boldsymbol{\beta},\boldsymbol{\eta}_0)-\rho(\boldsymbol{y}, \boldsymbol{X}\boldsymbol{\beta}_0,\boldsymbol{\eta}_0)]\mathbbm{1}_{\Lambda_{n}}(\vect{y},\vect{X})\}
	\label{A.6proof}
\end{align}

Using Cauchy-Schwarz inequality and A.2, we have that
\begin{align*}
\lvert \E\big[\rho(\boldsymbol{y}, \boldsymbol{X}\boldsymbol{\beta},\boldsymbol{\eta}_0)-\rho(\boldsymbol{y}, \boldsymbol{X}\boldsymbol{\beta}_0,\boldsymbol{\eta}_0) \mathbbm{1}_{\Lambda_{n}}(\vect{y},\vect{X})\big]\rvert \leq
    \E\big[\lvert \rho(\boldsymbol{y}, \boldsymbol{X}\boldsymbol{\beta},\boldsymbol{\eta}_0)-\rho(\boldsymbol{y}, \boldsymbol{X}\boldsymbol{\beta}_0,\boldsymbol{\eta}_0)\rvert \mathbbm{1}_{\Lambda_{n}}(\vect{y},\vect{X})\big] \\
	\leq k_n \E\big[\|\vect{X}\boldsymbol{1}_m\| \|\boldsymbol{\beta}-\boldsymbol{\beta}_0\|\big]\leq k_nm^{1/2}C^{1/2}b.
\end{align*}
Hence, we may chose a sufficiently large $D>0$ such that with $b$ chosen as in A.6, the second term in \eqref{A.6proof} can "catch up" with the first term in a way that makes the whole right hand side of \eqref{A.6proof} small enough leading to our desired results.
For Assumption A.7, we have that
\begin{align*}
	&\sum_{j=1}^m\Bigg[\E\Big\{\exp\big[(\vect{x}_{j}^T\boldsymbol{\beta}+\vect{z}_{j}^T\vect{u}+s_0)^2/2-(\vect{x}_{j}^T\boldsymbol{\beta}+\vect{z}_{j}^T\vect{u})^2/2\big]\\
	&~~~~~~~~~~~~~~~+\exp\big[(\vect{x}_{j}^T\boldsymbol{\beta}+\vect{z}_{j}^T\vect{u}-s_0)^2/2-(\vect{x}_{j}^T\boldsymbol{\beta}+\vect{z}_{j}^T\vect{u})^2/2\big]\Big\}\Bigg]\\
 &~~~~~~~~~~~~~~= \sum_{j=1}^m \E\Big[\exp\{s_0(\vect{x}_j^T\boldsymbol{\beta}+\vect{z}_j^T\vect{u})\}+\exp\{-s_0(\vect{x}_j^T\boldsymbol{\beta}+\vect{z}_j^T\vect{u})\}\Big]e^{s_0^2/2}\leq ms_3 e^{s_0^2/2},
\end{align*}
where the last inequality is due to B.2. Thus, by setting $s_1=ms_3e^{s_0^2/2}$, Assumption A.7 holds. This completes the proof of Lemma 3. $\hfill \square$

\section{Proof of Lemma 4}
For the first part of Condition A.3, we need to look at the Hessian of the loss function with respect to $\boldsymbol{\beta}$, which is given by differentiating the expression in brackets in \eqref{EQ:mdpdeforBetaEstimating};
\begin{align*}
 \nabla^2_{\boldsymbol{\beta}}\rho(\vect{y},\vect{X}\boldsymbol{\beta},\boldsymbol{\eta}) = -\big(\alpha+1\big)\frac{e^{-\frac{\alpha}{2}(\vect{y}-\vect{X}\boldsymbol{\beta})^T\vect{V}^{-1}(\vect{y}-\vect{X}\boldsymbol{\beta})}}{2\pi^{\frac{m\alpha}{2}}\lvert \vect{V}\rvert^{\frac{\alpha}{2}}} \Bigg[\alpha \vect{X}^T\vect{d}\vect{d}^T\vect{X}-\vect{X}^T\vect{V}^{-1}\vect{X}\Bigg],
\end{align*}
where $\vect{d}:= \vect{V}^{-1}(\vect{y}-\vect{X}\boldsymbol{\beta})$.
For the Hessian above to be positive semi-definite, the matrix inside the square brackets needs to be negative semi-definite, which holds due to Assumption B.3. For the second part of Condition 2, we need to find the generalized information matrix. 
Let $\text{vec}(\vect{V})$ be  the $m^2 \times 1$ column vector obtained by stacking the columns of the matrix $\vect{V}$ on top of one another. We then have the  gradient vector
\begin{align*}
    \nabla \rho = \begin{bmatrix}
\nabla_{\boldsymbol{\beta}}\rho\\
\nabla_{\boldsymbol{\eta}}\rho
\end{bmatrix}=
\begin{bmatrix}
2\cdot h(\boldsymbol{\beta},\boldsymbol{\eta})\cdot \vect{X}^T\vect{V}^{-1}(\vect{y}-\vect{X}\boldsymbol{\beta}) \\
        -\frac{\alpha}{2}L_1\lvert\vect{V}\rvert^{-\alpha/2}\text{vec}(\vect{V}^{-1})+h(\boldsymbol{\beta},\boldsymbol{\eta}) \cdot \big\{\text{vec}(\vect{V}^{-1}-\vect{V}^{-1}(\vect{y}-\vect{X}\boldsymbol{\beta})(\vect{y}-\vect{X}\boldsymbol{\beta})^T\vect{V}^{-1})\big\}
    \end{bmatrix}
\end{align*}
where $\nabla_{\boldsymbol{\beta}}\rho = \frac{\partial \rho(\vect{y},\vect{X}\boldsymbol{\beta},\boldsymbol{\eta})}{\partial {\boldsymbol{\beta}}}$ and  $\nabla_{\boldsymbol{\eta}} = \frac{\partial \rho\rho(\vect{y},\vect{X}\boldsymbol{\beta},\boldsymbol{\eta})}{\partial {\boldsymbol{\eta}}}$ and
\begin{align*}
h(\boldsymbol{\beta},\boldsymbol{\eta}) = 
    \frac{\alpha}{2}L_2\lvert\vect{V}\rvert^{-\alpha/2}\exp\big\{-\frac{\alpha}{2}(\vect{y}-\vect{X}\boldsymbol{\beta})^T\vect{V}^{-1}(\vect{y}-\vect{X}\boldsymbol{\beta})\big\}.
\end{align*}
The constants $L_1$ and $L_2$ are as defined in \eqref{EQ:MDPDElossfunctionLinearmixed}. The generaliz±ed information matrix is given by
\begin{align*}
\E(\nabla\rho\nabla\rho^T) = \E
    \begin{bmatrix}
\nabla_{\boldsymbol{\beta}}\rho\nabla_{\boldsymbol{\beta}}\rho^T & \nabla_{\boldsymbol{\beta}}\rho\nabla_{\boldsymbol{\eta}}\rho^T\\
\nabla_{\boldsymbol{\eta}}\rho\nabla_{\boldsymbol{\beta}}\rho^T & \nabla_{\boldsymbol{\eta}}\rho\nabla_{\boldsymbol{\eta}}\rho^T.
\end{bmatrix} 
\end{align*}

To get an expression for this matrix, we need the following results:
\begin{align*}
     \E\Big[\exp\big\{-\frac{\alpha}{2}(\vect{y}-\vect{X}\boldsymbol{\beta})^T\vect{V}^{-1}(\vect{y}-\vect{X}\boldsymbol{\beta})\big\}\vect{X}^T\vect{V}^{-1}(\vect{y}-\vect{X}\boldsymbol{\beta})\Big] \\= \int \frac{1}{(2\pi)^{m/2}}\frac{1}{\lvert \vect{V}\rvert^{1/2}}e^{-\frac{(1+\alpha)}{2}(\vect{y}-\vect{X}\boldsymbol{\beta})\vect{V}^{-1}(\vect{y}-\vect{X}\boldsymbol{\beta})}\vect{X}^T\vect{V}^{-1}(\vect{y}-\vect{X}\boldsymbol{\beta})d\vect{y} 
     \\= (1+\alpha)^{-m/2}\E[\vect{X}^T\vect{V}^{-1}(\widetilde{\vect{y}}-\vect{X}\boldsymbol{\beta})] = 0,
\end{align*}
where $\widetilde{\vect{y}} \sim \mathcal{N}(\vect{X}\boldsymbol{\beta},(1+\alpha)^{-1}\vect{V})$. 
With same arguments, it can be shown that  
\begin{align*}
     \E\Big[\exp\big\{-\alpha(\vect{y}-\vect{X}\boldsymbol{\beta})^T\vect{V}^{-1}(\vect{y}-\vect{X}\boldsymbol{\beta})\big\}\vect{X}^T\vect{V}^{-1}(\vect{y}-\vect{X}\boldsymbol{\beta})\Big] = 0.
\end{align*}
Now, 
\begin{align*}
    \E\Big[\exp\big\{-\alpha(\vect{y}-\vect{X}\boldsymbol{\beta})^T\vect{V}^{-1}(\vect{y}-\vect{X}\boldsymbol{\beta})\big\}\vect{X}^T\vect{V}^{-1}(\vect{y}-\vect{X}\boldsymbol{\beta})\cdot\text{vec}(\vect{V}^{-1}(\vect{y}-\vect{X}\boldsymbol{\beta})(\vect{y}-\vect{X}\boldsymbol{\beta})^T\vect{V}^{-1})^T\Big] \\
    = \E\Big[\exp\big\{-\alpha(\vect{y}-\vect{X}\boldsymbol{\beta})^T\vect{V}^{-1}(\vect{y}-\vect{X}\boldsymbol{\beta})\big\}\vect{X}^T\vect{V}^{-1}(\vect{y}-\vect{X}\boldsymbol{\beta})
    (\vect{y}-\vect{X}\boldsymbol{\beta})^T\otimes (\vect{y}-\vect{X}\boldsymbol{\beta})^T
\vect{V}^{-1}\otimes\vect{V}^{-1}\Big] \\
= \vect{X}^T\vect{V}^{-1}\E\Big[\exp\big\{-\alpha(\vect{y}-\vect{X}\boldsymbol{\beta})^T\vect{V}^{-1}(\vect{y}-\vect{X}\boldsymbol{\beta})\big\}(\vect{y}-\vect{X}\boldsymbol{\beta})
    (\vect{y}-\vect{X}\boldsymbol{\beta})^T\otimes (\vect{y}-\vect{X}\boldsymbol{\beta})^T\Big]
\vect{V}^{-1}\otimes\vect{V}^{-1}\\
    = (1+2\alpha)^{-m/2}\vect{X}^T\vect{V}^{-1}\E[(\check{\vect{y}}-\vect{X}\boldsymbol{\beta})
    (\check{\vect{y}}-\vect{X}\boldsymbol{\beta})^T\otimes (\check{\vect{y}}-\vect{X}\boldsymbol{\beta})^T]\vect{V}^{-1}\otimes\vect{V}^{-1} = 0,
\end{align*}
where $\otimes$ denotes the Kronecker product and $\check{\vect{y}}\sim \mathcal{N}(\vect{X}\boldsymbol{\beta},(1+2\alpha)^{-1}\vect{V})$. The last equality is because the expression inside the square brackets is the third central moment of a multivariate normal distribution which equals $0$. 
Hence, 
\begin{align*}
\E[\nabla_{\boldsymbol{\beta}}\rho\nabla_{\boldsymbol{\eta}}\rho^T]= \E[\nabla_{\boldsymbol{\eta}}\rho\nabla_{\boldsymbol{\beta}}\rho^T] = 0.
\end{align*}

Furthermore, we have that 
\begin{align*}
\E[\nabla_{\boldsymbol{\beta}}\rho\nabla_{\boldsymbol{\beta}}\rho^T]= \alpha^2L_2^2\lvert\vect{V}\rvert^{-\alpha}\vect{X}^T\vect{V}^{-1}\E\big[\exp\big\{-\alpha(\vect{y}-\vect{X}\boldsymbol{\beta})^T\vect{V}^{-1}(\vect{y}-\vect{X}\boldsymbol{\beta})\big\}(\vect{y}-\vect{X}\boldsymbol{\beta})(\vect{y}-\vect{X}\boldsymbol{\beta})^T\big]\vect{V}^{-1}\vect{X},
\end{align*}
where
\begin{align*}
E\big[\exp\big\{-\alpha(\vect{y}-\vect{X}\boldsymbol{\beta})^T\vect{V}^{-1}(\vect{y}-\vect{X}\boldsymbol{\beta})\big\}(\vect{y}-\vect{X}\boldsymbol{\beta})(\vect{y}-\vect{X}\boldsymbol{\beta})^T\big] \\= (1+2\alpha)^{-m/2}\E[(\check{\vect{y}}-\vect{X}\boldsymbol{\beta})(\check{\vect{y}}-\vect{X}\boldsymbol{\beta})^T] \\
= (1+2\alpha)^{-(1+m/2)}\vect{V}.
\end{align*}
Hence,
\begin{align*}
\E[\nabla_{\boldsymbol{\beta}}\rho\nabla_{\boldsymbol{\beta}}\rho^T] = \alpha^2L_2^2\lvert\vect{V}\rvert^{-\alpha} (1+2\alpha)^{-(1+m/2)}\vect{X}^T\vect{V}^{-1}\vect{X} \succ \vect{0},
\end{align*}
because $\vect{V}$ is positive definite by B.1. 

Finally, for $\E[\nabla_{\boldsymbol{\eta}}\rho\nabla_{\boldsymbol{\eta}}\rho^T]$, we need the following results:
\begin{align*}
       \E\Big[\exp\big\{-\frac{\alpha}{2}(\vect{y}-\vect{X}\boldsymbol{\beta})^T\vect{V}^{-1}(\vect{y}-\vect{X}\boldsymbol{\beta})\big\}\Big] = (1+\alpha)^{-m/2},
\end{align*}
because it is the moment-generating function of a Chi-squared distributed variable with $m$ degrees of freedom. Furthermore, we have that
\begin{align*}
    \E\Big[\exp\big\{-\frac{\alpha}{2}(\vect{y}-\vect{X}\boldsymbol{\beta})^T\vect{V}^{-1}(\vect{y}-\vect{X}\boldsymbol{\beta})\big\}\text{vec}(\vect{V}^{-1}(\vect{y}-\vect{X}\boldsymbol{\beta})(\vect{y}-\vect{X}\boldsymbol{\beta})^T\vect{V}^{-1})\Big]\\
    = \E\Big[\exp\big\{-\frac{\alpha}{2}(\vect{y}-\vect{X}\boldsymbol{\beta})^T\vect{V}^{-1}(\vect{y}-\vect{X}\boldsymbol{\beta})\big\}\vect{V}^{-1}\otimes\vect{V}^{-1}\text{vec}\big((\vect{y}-\vect{X}\boldsymbol{\beta})(\vect{y}-\vect{X}\boldsymbol{\beta})^T\big)
\Big]\\
=  \vect{V}^{-1}\otimes\vect{V}^{-1}\text{vec}\Big(\E\big[\widetilde{\vect{y}}-\vect{X}\boldsymbol{\beta})(\widetilde{\vect{y}}-\vect{X}\boldsymbol{\beta})^T\big]\Big)\\
= (1+\alpha)^{-1}\vect{V}^{-1}\otimes\vect{V}^{-1}\text{vec}(\vect{V}).
\end{align*}
Similarly,
\begin{align*}
      \E\Big[\exp\big\{-\alpha(\vect{y}-\vect{X}\boldsymbol{\beta})^T\vect{V}^{-1}(\vect{y}-\vect{X}\boldsymbol{\beta})\big\}\text{vec}(\vect{V}^{-1}(\vect{y}-\vect{X}\boldsymbol{\beta})(\vect{y}-\vect{X}\boldsymbol{\beta})^T\vect{V}^{-1})\Big]\\
      = (1+2\alpha)^{-1}\vect{V}^{-1}\otimes\vect{V}^{-1}\text{vec}(\vect{V}) \\
\end{align*}
Lastly, define $\vect{d} = \vect{V}^{-1}({\vect{y}}-\vect{X}\boldsymbol{\beta})$. Then we also need to find
\begin{align*}
    \E\Big[\exp\big\{-\alpha ({\vect{y}}-\vect{X}\boldsymbol{\beta})^T\vect{V}^{-1}({\vect{y}}-\vect{X}\boldsymbol{\beta})\big\}\text{vec}(\vect{d}\vect{d}^T)\text{vec}(\vect{d}\vect{d}^T)^T\Big]
    =\E\Big[e^{-\alpha ({\vect{y}}-\vect{X}\boldsymbol{\beta})^T\vect{V}^{-1}({\vect{y}}-\vect{X}\boldsymbol{\beta})}\vect{d}\vect{d}^T\otimes\vect{d}\vect{d}^T\Big],
\end{align*}
which is positive definite by Assumption B.5. Combining all of the results above and B.4, we get that each term in $\E[\nabla_{\boldsymbol{\eta}}\rho\nabla_{\boldsymbol{\eta}}\rho^T]$ is positive semi-definite, and hence the sum is positive semi-definite, i.e., Assumption A.2. is fulfilled given B.3--B.5. For Condition A.4, we may use that if A.3 holds, then there is a neighborhood $\boldsymbol{\mathcal{B}}$ around $(\boldsymbol{\beta}_0,\boldsymbol{\eta}_0)$ in which the expected value of the loss function is strongly convex with a positive constant $V_m$, because $\E[\nabla\rho(\vect{y},\vect{X}\boldsymbol{\beta}_0,\boldsymbol{\eta}_0)] = 0$. Similarly to MLE, condition A.5 holds since $\rho(\vect{y},\vect{X}\boldsymbol{\beta},\boldsymbol{\eta})$ is continuously differentiable for all $\vect{X}\boldsymbol{\beta}$, $\boldsymbol{\beta} \in \boldsymbol{\mathcal{B}}$. Finally, Assumption A.6 holds for the same reasons as MLE.

\hfill$\square$

\end{document}